\documentclass[12pt, a4paper,reqno]{amsart}

\usepackage{amsmath}
\usepackage{amsfonts}
\usepackage{amssymb}
\usepackage{amsthm}
\usepackage[active]{srcltx}

\numberwithin{equation}{section}

\newtheorem{theorem}[equation]{Theorem}
\newtheorem{lemma}[equation]{Lemma}

\theoremstyle{definition}
\newtheorem{definition}[equation]{Definition}

\theoremstyle{remark}

\def\kint_#1{\mathchoice%
          {\mathop{\kern 0.2em\vrule width 0.6em height 0.69678ex depth -0.58065ex
                  \kern -0.8em \intop}\nolimits_{\kern -0.4em#1}}%
          {\mathop{\kern 0.1em\vrule width 0.5em height 0.69678ex depth -0.60387ex
                  \kern -0.6em \intop}\nolimits_{#1}}%
          {\mathop{\kern 0.1em\vrule width 0.5em height 0.69678ex depth -0.60387ex
                  \kern -0.6em \intop}\nolimits_{#1}}%
          {\mathop{\kern 0.1em\vrule width 0.5em height 0.69678ex depth -0.60387ex
                  \kern -0.6em \intop}\nolimits_{#1}}}
\def\vintslides_#1{\mathchoice%
          {\mathop{\kern 0.1em\vrule width 0.5em height 0.697ex depth -0.581ex
                  \kern -0.6em \intop}\nolimits_{\kern -0.4em#1}}%
          {\mathop{\kern 0.1em\vrule width 0.3em height 0.697ex depth -0.604ex
                  \kern -0.4em \intop}\nolimits_{#1}}%
          {\mathop{\kern 0.1em\vrule width 0.3em height 0.697ex depth -0.604ex
                  \kern -0.4em \intop}\nolimits_{#1}}%
          {\mathop{\kern 0.1em\vrule width 0.3em height 0.697ex depth -0.604ex
                  \kern -0.4em \intop}\nolimits_{#1}}}

\newcommand{\R}{\mathbb{R}}
\newcommand{\Rn}{\mathbb{R}^d}

\newcommand{\esssup}{\operatornamewithlimits{ess\, sup}}
\newcommand{\essinf}{\operatornamewithlimits{ess\,inf}}
\newcommand{\essosc}{\operatornamewithlimits{ess\,osc}}

\renewcommand{\div}{\nabla \cdot}

\renewcommand{\l}{\left}
\renewcommand{\r}{\right}

\def\Xint#1{\mathchoice
{\XXint\displaystyle\textstyle{#1}}%
{\XXint\textstyle\scriptstyle{#1}}%
{\XXint\scriptstyle\scriptscriptstyle{#1}}%
{\XXint\scriptscriptstyle\scriptscriptstyle{#1}}%
\!\int}

\def\XXint#1#2#3{{\setbox0=\hbox{$#1{#2#3}{\int}$}
\vcenter{\hbox{$#2#3$}}\kern-.5\wd0}}

\def\dashint{\Xint-}

\title[H\"older regularity]{A note on the proof of H\"older continuity to weak solutions of elliptic equations}
\author{Juhana Siljander}
\address{Aalto University, Institute of Mathematics, P.O. Box 11100, FI-00076 Aalto, Finland.}
\email{juhana.siljander@tkk.fi}

\begin{document}

\begin{abstract}
By borrowing ideas from the parabolic theory, we use a combination of De Giorgi's and Moser's methods to give some remarks on the proof of H\"older continuity of weak solutions of elliptic equations.
\end{abstract}
%\date{\today}


%\subjclass[2010]{Primary 35J92. Secondary 35B65}

\maketitle

\section{Introduction}
We will present some observations on the proof of H\"older continuity of weak solutions to equations of type
\begin{equation}\label{equation}
\div{\mathcal{A}(x,u,\nabla u)}=0.
\end{equation}
This kind of elliptic equations are, of course, well-studied and there are many beautiful arguments for the H\"older regularity of their solutions. As it is well known, the problem was first solved independently by Ennio De Giorgi~\cite{DeGi57} and John Nash~\cite{Nash58}. After J\"urgen Moser used his iteration method for proving the supremum estimate, methods based on Harnack's inequalities were found as well~\cite{Mose71},~\cite{Trud68},~\cite{Mose64},~\cite{Mose64corr},~\cite{Mose60}.

Although the elliptic case is very well understood nowadays, the parabolic case seems to be more involved. In particular, there seems to be only one method for proving the continuity result for parabolic equations~\cite{DiBe93}. Consequently, a lot of research has been done for understanding the parabolic theory.

Using the ideas developed for the parabolic equations, we will give some remarks on the proof of H\"older continuity in the elliptic case. More presicely, we combine the De Giorgi method, in a form used in the parabolic setting, with Moser's iteration and a crossover lemma to give a proof for the regularity theorem.

The argument is formulated for a general Borel measure which is assumed to satisfy the doubling condition and to support a weak Poincar\'e inequality. These together are known to imply a Sobolev inequality which is the crucial tool we use. Regularity arguments for elliptic equations in the weighted case have been studied, for instance, by Fabes, Kenig and Serapioni in~\cite{FabeKeniSera82}. For further aspects of the theory see the classical book by Ladyzhenskaya and Uraltseva~\cite{LadyUral68}.

\section{Preliminaries} \label{section:preliminaries}

Let $\mu$ be a Borel measure and $\Omega$ an open set in $\Rn$. The Sobolev space $H^{1,p}(\Omega)$ is defined to be the completion of $C^\infty(\Omega)$ with respect to the Sobolev norm
\[
\|u\|_{1,p,\Omega}=\left(\int_\Omega |u|^p + |\nabla u|^p \, d\mu \right)^{1/p}.
\]

A function $u$ belongs to the local Sobolev space $H_{loc}^{1,p}(\Omega)$ if it belongs to $H^{1,p}(\Omega')$ for every $\Omega' \Subset \Omega$. Moreover, the Sobolev space with zero boundary values is defined as the completion of $C_0^\infty(\Omega)$ with respect to the Sobolev norm. For more properties of Sobolev spaces, see e.g. \cite{HeinKilpMart93} or \cite{AdamFour03}.

Assume that $\mathcal{A}:\Rn\times\R\times\Rn\rightarrow \Rn$ is a function such that $\mathcal{A}(\cdot,\zeta,\xi)$ is measurable for every $(\zeta,\xi)\in \R\times \Rn$ and $\mathcal{A}(x,\cdot,\cdot)$ is continuous for all $x\in \Omega$. Suppose also that for some $\mathcal{A}_0\ge 0$ and $\mathcal{C}_0>0$ we have
\[
 |\mathcal{A}(x,\zeta,\xi)|\le\mathcal{A}_0|\xi|^{p-1}
\]
and
\[
 \mathcal{A}(x,\zeta,\xi)\cdot\xi\ge \mathcal{C}_0|\xi|^{p}.
\]
A weak solution for equation~\eqref{equation} is defined as follows.

\begin{definition}
  A function $u \in H_{loc}^{1,p}(\Omega)$ is a
  weak solution of equation \eqref{equation} in $\Omega$ if it satisfies the integral equality
\begin{equation}\label{weak_solution}
\begin{split}
\int_{\Omega} \mathcal{A}(x,u,\nabla u)\cdot \nabla \phi
\, d\mu \ =& \ \ 0
\end{split}
\end{equation}
for all $\phi \in C_0^\infty(\Omega)$. If the equality in this definition is replaced by $\ge$ ($\le$) and the inequality holds for every nonnegative $\phi \in C_0^\infty(\Omega)$ we say that the function is a supersolution (subsolution).
\end{definition}

The measure $\mu$ is said to be doubling if there is a universal constant $D_0\ge 1$ such that
\[
 \mu(B(x,2r))\le D_0 \mu(B(x,r))
\]
for all $B(x,2r) \subset\Omega$. Here $B(x,r)$ denotes the standard open ball in $\Rn$
\[
B(x,r)=\{y\in \Rn : |y-x|<r\}.
\]
We will also use the notation
\[
 B(r):=B(0,r).
\]

The dimension related to the doubling measure is defined by $d_\mu:=\log_2 D_0$. Note that in the case of the Lebesgue measure $d_\mathcal{L}=d$. The measure is said to support a weak $(1,p)$-Poincar\'e inequality if there exist constants $P_0>0$ and $\tau\ge 1$ such that
\begin{equation}\label{poincare}
\dashint_{B(x,r)}|u-u_{B(x,r)}| \, d\mu \le P_0 r\left(\dashint_{B(x,\tau r)} |\nabla u|^p \, d\mu\right)^{1/p}
\end{equation}
for every $u\in H^{1,p}(\Omega)$ and $B(x,\tau r)\subset \Omega$. Here we used the notation
\[
 u_{B(x,r)}=\dashint_{B(x,r)} u \, d\mu = \frac{1}{\mu(B(x,r))}\int_{B(x,r)} u \, d\mu.
\]
The word weak refers to the constant $\tau\ge 1$. If the inequality \eqref{poincare} is true for $\tau=1$ we say that the measure supports a $(1,p)$-Poincar\'e inequality.

%From now on will assume the measure to be doubling and to support the weak $(1,p)$-Poincar\'e inequality. Moreover, we assume the measure to be non-trivial in the sense that the measure of every non-empty open set is strictly positive and the measure of every bounded set is finite.

It is known that the weak $(1,p)$-Poincar\'e inequality and the doubling condition imply a Sobolev embedding.

\begin{theorem}\label{Sobolev}
Suppose $u \in H_0^{1,p}(B(x,r))$. Then there is a constant $C>0$ such that
\[
\left(\dashint_{B(x,r)} |u|^\kappa \, d\mu\right)^{1/\kappa}\le Cr\left(\dashint_{B(x,r)} |\nabla u|^p\, d\mu\right)^{1/p}
\]
where 
\[
\kappa=\begin{cases}\frac{d_\mu p}{d_\mu-p}, \quad \text{for} \quad 1<p<d_\mu \\ 2p, \quad \text{otherwise}. \end{cases}
\]
\begin{proof} See for example \cite{KinnShan01}.
\end{proof}
\end{theorem}

We will also need the following lemma.

\begin{lemma}\label{geometric_convergence}
Let $\{Y_n\}, n=0,1,2,\dots$, be a sequence of positive numbers, satisfying
\[
Y_{n+1}\le Cb^nY_n^{1+\alpha}
\]
where $C,b>1$ and $\alpha>0$. Then $\{Y_n\}$ converges to zero as $n\rightarrow\infty$, provided
\[
Y_0\le C^{-1/\alpha}b^{1-\alpha^2}.
\]
\begin{proof} For the proof we refer to \cite{DiBe93}.
\end{proof}
\end{lemma}

Our main theorem is the following well-known regularity result. The observations we make lie in the proof of the claim. More precisely, to deduce the claim we use a combination of De Giorgi's method and Moser's iteration scheme together with Chebyshev's inequality.

\begin{theorem}
Suppose $\mu$ is a doubling measure which supports a weak $(1,p)$-Poincar\'e inequality. Let $u\in H_{loc}^{1,p}(\Omega)$ be a weak solution of equation~\eqref{equation}. Then $u$ is locally H\"older continuous.
\end{theorem}

We will prove the H\"older continuity of the solution in a neighborhood of an arbitrary point. Since the equation is translation invariant, for simplicity of notation, we can assume this point to be the origin.

\pagebreak

\section{Estimates for weak solutions}

Let us start by stating some classical lemmata.

\begin{lemma}[Caccioppoli]\label{energy}
Let $u\ge 0$ be a weak subsolution for equation~\eqref{equation} in $\Omega$. Then there exists a constant $C=C(p,\mathcal{A}_0,\mathcal{C}_0)>0$ such that for every $k\ge 0$ and $\varphi\in C_0^\infty(\Omega)$ we have
\[
\int_\Omega |\nabla(u-k)_+|^p\varphi^p \, d\mu \le C\int_\Omega (u-k)_+^p|\nabla\varphi|^p \, d\mu.
\]

\begin{proof}
The result follows by choosing the test function $\phi=(u-k)_+\varphi^p$ in the definition of a weak solution. For details see~\cite{MalyZiem97}.
\end{proof}
\end{lemma}

\begin{lemma}[Crossover]\label{crossover}
Let $u\ge 0$ be a weak supersolution for equation~\eqref{equation} in $\Omega$ and let $B(r)\Subset\Omega$. Then there exist constants $C$ and $\delta>0$ such that
\[
\l(\dashint_{B(r)} u^{-\delta}\, d\mu\r)^{1/\delta}\le C \l(\dashint_{B(r)} u^{\delta}\, d\mu\r)^{-1/\delta}.
\]
\begin{proof}
For the proof we refer to~\cite{HeinKilpMart93}.
\end{proof}
\end{lemma}

\begin{lemma}\label{esssup}
 Let $u \ge 0$ be a weak solution of equation~\eqref{equation} in $\Omega$ and let $B(r) \Subset \Omega$. Then for every $\delta>0$ there exists a constant $C>0$ such that 
\[
 \esssup_{B(r/2)}{u}\le C \l(\dashint_{B(r)} u^\delta \, d\mu\r)^{1/\delta}.
\]
\begin{proof}
 The result follows by standard iteration techniques, see~\cite{HeinKilpMart93}.
\end{proof}

\end{lemma}

\section{H\"older Continuity}
%\begin{lemma}
%Let $u\ge 0$ be a weak solution for equation~\eqref{equation}. Then $v:=1/u$ is a subsolution of an equation of similar type. Moreover, $v$ is locally bounded.
%\begin{proof}
%\end{proof}
%\end{lemma}

Let $r>0$ and denote
\[
r_n:=\frac{r}{2}+\frac{r}{2^{n+1}},\qquad B_{n}:=B(r_n)
\]
and
\[
A_{n}:=\left\{x\in B_{n}: u(x)>k_n\right\}
\]
where

%epsilon-version
%\[
% k_n=\esssup_{B(r)}{u} - \frac{\essosc_{B(r)}{u}}{2^{\lambda+1}}-\frac{\essosc_{B(r)}{u}}{2^{\lambda+n+1}}+\eps
%\]
%non-epsilon version
\[
 k_n:=\esssup_{B(r)}{u} - \frac{\essosc_{B(r/2)}{u}}{2^{\lambda+1}}-\frac{\essosc_{B(r/2)}{u}}{2^{\lambda+n+1}}
\]
for $n=0,1,2,\dots$.% and $\eps>0$. 

\begin{lemma}\label{main_lemma}
Let $u$ be a weak subsolution of equation~\eqref{equation} in B(r). Then there exists a constant $C>0$ such that
\[
\frac{\mu(A_{n+1})}{\mu(B_{n+1})}\le C 4^{n\kappa(1+1/p)}\l(\frac{\mu(A_n)}{\mu(B_n)}\r)^{\kappa/p}.
\]
\begin{proof}
Choose the cut off function $\varphi_n\in C_0^\infty(B_n)$ such that $\varphi_n=1$ in $B_{n+1}$ and
\[
|\nabla \varphi_n|\le \frac{C2^n}{r}, \quad n=1,2,\dots.
\]

Using the doubling property of the measure together with Sobolev's inequality (Theorem~\ref{Sobolev}) and the Caccioppoli inequality (Lemma~\ref{energy}) gives
\begin{align*}
&\dashint_{B_{n+1}} (u-k_n)_+^{\kappa}\, d\mu  \\
&\le \dashint_{B_{n+1}} (u-k_n)_+^\kappa\varphi_n^\kappa\, d\mu\le \frac{\mu(B_{n})}{\mu(B_{n+1})}\dashint_{B_{n}} (u-k_n)_+^\kappa\varphi_n^\kappa\, d\mu \\
&\le Cr^\kappa\l(\dashint_{B_n} |\nabla (u-k_n)_+\varphi_n|^{p}\, d\mu\r)^{\kappa/p} \\
&\le Cr^\kappa \l(\dashint_{B_{n}} |\nabla (u-k_n)_+|
^p\varphi_n^p+(u-k_n)_+^p|\nabla \varphi_n|^p\, d\mu\r)^{\kappa/p} \\
&\le Cr^\kappa \l(\dashint_{B_{n}} (u-k_n)_+^p|\nabla \varphi_n|^p\, d\mu\r)^{\kappa/p} \\
&\le C 2^{n\kappa/p}\l(\dashint_{B_{n}} (u-k_n)_+^p\, d\mu\r)^{\kappa/p} \\
%non-epsilon
&\le C 2^{n\kappa/p}\l(\frac{\essosc_{B(r/2)}{u}}{2^{\lambda}}\r)^\kappa \l(\frac{\mu(A_n)}{\mu(B_n)}\r)^{\kappa/p}.
%epsilon
%&\le C 2^{n\kappa/p}\l(\frac{\essosc_{B(r)}{u}}{2^{\lambda}}-\eps\r)^\kappa \l(\frac{\mu(A_n)}{\mu(B_n)}\r)^{\kappa/p}.
\end{align*}
On the other hand,
\[
\dashint_{B_{n+1}} (u-k_n)_+^{\kappa}\, d\mu \ge \frac{\mu(A_{n+1})}{\mu(B_{n+1})}\l(\frac{\essosc_{B(r/2)}{u}}{2^{\lambda+n+2}}\r)^\kappa.
\]
These together give
\[
\frac{\mu(A_{n+1})}{\mu(B_{n+1})}\le C 4^{n\kappa(1+1/p)}\l(\frac{\mu(A_n)}{\mu(B_n)}\r)^{\kappa/p},
\]
as required.
\end{proof}
\end{lemma}

Now by Lemma~\ref{geometric_convergence} we have $\mu(A_n)/\mu(B_n)\rightarrow 0$ as $\rightarrow 0$, provided 
\begin{equation}\label{smallness}
\frac{\mu(A_0)}{\mu(B_0)}\le C^{-1/(\kappa/p-1)} 4^{\kappa(1-(1-\kappa/p)^2)}.
\end{equation}
Next we turn to prove that this will, indeed, be satisfied for some suitably chosen $\lambda>0$.

\begin{lemma}
Let $u$ be a weak solution of equation~\eqref{equation} in $B(3r)$. Then there exists a constant $\lambda_0:=\lambda>0$ such that~\eqref{smallness} holds. Recall that $A_0$ depends on $\lambda$.

\begin{proof}
Now by Chebyshev's inequality we have

%epsilon

%\begin{align*}
%\mu(A_0)&=\mu(\{x\in B(r): u>\esssup_{B(r)}{u} - \frac{\essosc_{B(r)}{u}}{2^{\lambda}}+\eps\}) \\
%&=\mu(\{x\in B(r): \frac{\essosc_{B(r)}{u}}{2^{\lambda}}>\esssup_{B(r)}{u}-u +\eps \}) \\
%&\le \frac{\essosc_{B(r)}{u}}{2^{\lambda}}\int_{A_0}\frac{1}{\esssup_{B(r)}{u}-u+\eps}.
%\end{align*}

%Since 
%\[
%\esssup_{B(r)}{u}-u+\eps\ge\eps> 0
%\]
%is a weak solution for equation~\eqref{equation} in $B(r)$, 
%\[
%\frac{1}{\esssup_{B(r)}{u}-u+\eps}
%\]
%is a subsolution for an equation of similar type. Subsolutions are locally bounded so we have
%\[
%\int_{A_0}\frac{1}{\esssup_{B(r)}{u}-u+\eps}\le C_B
%\]
%for some constant $C_B$ which depends only on the ball $B(r)\subset \Omega$. In particular $C_B$ will not depend on $\eps$! Now choosing $\lambda$ large enough guarantees~\eqref{smallness}. This gives that for some $\lambda_0>0$ which depends only upon the data and the ball $B(r)\subset \Omega$ we are in, we have
%\[
%u\le \esssup_{B(r)}{u} - \frac{\essosc_{B(r)}{u}}{2^{\lambda_0}}+\eps\quad\text{a.e in}\quad B(r/2).
%\]
%This is uniform for every $\eps>0$ so we can let $\eps\rightarrow 0$ and hence we will get 
%\[
%u\le \esssup_{B(r)}{u} - \frac{\essosc_{B(r)}{u}}{2^{\lambda_0}}\quad\text{a.e in}\quad B(r/2).
%\]
%H\"older continuity follows by a standard argument.

%\end{proof}
%\end{theorem}

%non-epsilon

\begin{align*}
\frac{\mu(A_0)}{\mu(B_0)}&=\mu(\{x\in B(r): u>\esssup_{B(r)}{u} - \frac{\essosc_{B(r/2)}{u}}{2^{\lambda}}\})/\mu(B_0) \\
&=\mu(\{x\in B(r): \frac{\essosc_{B(r/2)}{u}}{2^{\lambda}}>\esssup_{B(r)}{u}-u  \})/\mu(B_0)\\
&\le \l(\frac{\essosc_{B(r/2)}{u}}{2^{\lambda}}\r)^\delta\dashint_{B(r)}\l(\frac{1}{\esssup_{B(r)}{u}-u}\r)^\delta \, d\mu
\end{align*}
where $\delta>0$ is to be determined shortly. Since 
\[
\esssup_{B(r)}{u}-u\ge 0
\]
is a weak solution of equation~\eqref{equation} in $B(3r)$, by the Crossover lemma (Lemma~\ref{crossover}) we have
\[
\dashint_{B(r)}\l(\frac{1}{\esssup_{B(r)}{u}-u}\r)^{\delta} \, d\mu \le C \l(\dashint_{B(r)} (\esssup_{B(r)}{u}-u)^{\delta} \, d\mu\r)^{-1}
\]
for all small enough $\delta>0$. By Lemma~\ref{esssup} we obtain
\begin{align*}
 \l(\dashint_{B(r)} (\esssup_{B(r)}{u}-u)^{\delta} \, d\mu\r)^{1/\delta}&\ge\frac{1}{C} \esssup_{B(r/2)}{(\esssup_{B(r)}{u}-u)} \\
&=\frac{1}{C}\l(\esssup_{B(r)}{u}-\essinf_{B(r/2)}{u}\r) \\
&\ge \frac{1}{C}\essosc_{B(r/2)}{u}.
\end{align*}
Consequently,
\[
 \frac{\mu(A_0)}{\mu(B_0)}\le \frac{C}{2^{\delta\lambda}}.
\]
Choosing $\lambda$ large enough finishes the proof.

% where the right hand side is clearly bounded, provided $u$ is not constant in $B(r)$. Moreover, the bound is uniform inside $B(r)$ in the sense that for any compact set $K\subset B(r)$ the bound can be chosen to be the same. Now choosing $\lambda$ large enough guarantees~\eqref{smallness} uniformly in $B(r)$. 
%\begin{align*}
%&\dashint_{B_{n+1}} (u-k_n)_+^{p\gamma}\, d\mu  \\
%&\le \l(\dashint_{B_{n+1}} (u-k_n)_+^p\varphi^p\, d\mu\r)^{(\kappa-p)/\kappa}\l(\dashint_{B_{n+1}} (u-k_n)_+^\kappa\varphi^\kappa\, d\mu\r)^{p/\kappa} \\
%&\le C\l(\dashint_{B_{n}} (u-k_n)_+^p\varphi^p\, d\mu\r)^{(\kappa-p)/\kappa}\l(\dashint_{B_{n}} (u-k_n)_+^\kappa\varphi^\kappa\, d\mu\r)^{p/\kappa} \\
%&\le Cr^p\l(\dashint_{B_{n}} (u-k_n)_+^p\varphi^p\, d\mu\r)^{(\kappa-p)/\kappa}\dashint_{B_n} |\nabla (u-k_n)_+\varphi|^{p}\, d\mu \\
%&\le Cr^p \l(\dashint_{B_{n}} (u-k_n)_+^p\varphi^p\, d\mu+|\nabla (u-k_n)_+|
%^p\varphi^p+(u-k_n)_+^p|\nabla \varphi|^p\, d\mu\r)^{\gamma} \\
%&\le Cr^p \l(\dashint_{B_{n}} (u-k_n)_+^p\varphi^p\, d\mu+(u-k_n)_+^p|\nabla \varphi|^p\, d\mu\r)^{\gamma}
%\end{align*}

\end{proof}
\end{lemma}

%\begin{remark}
%By Lemma~\ref{crossover} we also have that if $u\ge 0$ is $\mathcal{A}$-harmonic and not identically zero in $\Omega$, then $\mu(x\in \Omega: u=0)=0$. Moreover, since $f(t)=1/t$ is convex when $t>0$, we can deduce that $1/u$ is $\mathcal{A}$-subharmonic in the set $\{u>0\}$.
%\end{remark}

Now the H\"older estimate follows from the previous result by standard measures. For the sake of completeness we recall the argument in the form of the following theorem.
\begin{theorem}
Let $u$ be a weak solution of equation~\eqref{equation} in $B(3r)$ and let $x,y \in B(r)$ and $r>|x-y|/2$. Then there exist constants $C>0$ and $0<\alpha<1$ such that
\[
|u(x)-u(y)| \le C\l(\frac{|x-y|}{r}\r)^\alpha\esssup_{B(3r)}{|u|}.
\]

\begin{proof}
% By the previous Lemma for some $\lambda_0>0$ %which depends only upon the data and the ball $B(r)\subset \Omega$ we are in, we have
% we have
% \begin{equation}\label{decrease}
% u\le \esssup_{B(r)}{u} - \frac{\essosc_{B(r/2)}{u}}{2^{\lambda_0}}\quad\text{a.e in}\quad B(r/2).
% \end{equation}
Either
\[
 \essosc_{B(r/2)}{u}\le\frac{1}{2}\essosc_{B(r)}{u}
\]
or by using the previous Lemma together with Lemmas~\ref{main_lemma} and~\ref{geometric_convergence} we obtain
\begin{equation}\label{decrease}
u\le \esssup_{B(r)}{u} - \frac{\essosc_{B(r)}{u}}{2^{\lambda_0+2}}\quad\text{a.e in}\quad B(r/2)
\end{equation}
for some $\lambda_0>0$, which only depends on the data. Now by subtracting
\(
\essinf_{B(r/2)}{u}
\)
from both sides of~\eqref{decrease} we obtain
\begin{equation}\label{reduction}
\essosc_{B(r/2)}{u}\le \l(1-\frac{1}{2^{\lambda_0+2}}\r)\essosc_{B(r)}{u}.
\end{equation}
We conclude that in any case~\eqref{reduction} is true.

Let $\gamma=(1-1/2^{\lambda_0+2})$, $0<r<R$ and choose $i$ such that
\begin{equation}\label{i}
\frac{R}{2^{i+1}}\le r \le \frac{R}{2^{i}}.
\end{equation}
 Now this together with an iteration of~\eqref{reduction} gives
\begin{align*}
\essosc_{B(r)}{u}\le \essosc_{B(R/2^i)}{u} \le \gamma^i\essosc_{B(R)}{u} \le C \l(\frac{r}{R}\r)^\alpha\essosc_{B(R)}{u}
\end{align*}
where
\[
\alpha=-\frac{\log{\gamma}}{\log{2}}.
\]
Let now $x,y \in B(r)$ and, further, let $R=2r>|x-y|$. Now we have
\begin{align*}
|u(x)-u(y)|&\le \essosc_{B((x+y)/2,|x-y|)}{u}\\
&\le C\l(\frac{|x-y|}{R}\r)^\alpha\essosc_{B((x+y)/2,R)}{u} \\
&\le C\l(\frac{|x-y|}{r}\r)^\alpha\esssup_{B(3r)}{|u|},
\end{align*}
as required.
\end{proof}
\end{theorem}

\bibliography{citations}
\bibliographystyle{plain}

\end{document}